\documentclass[12pt]{amsart}
\usepackage{amsfonts}
\usepackage{amssymb}
\usepackage{amsmath}
\usepackage[all]{xy}
\usepackage[english]{babel}
\usepackage{xspace}
\usepackage{times}
\usepackage{amsthm}
\usepackage{amscd}
\usepackage{amsbsy}
\usepackage{mathrsfs}
\usepackage{relsize}
\usepackage[center,small,sc]{caption}
\usepackage{pgf}
\usepackage{lmodern}
\usepackage{xcolor}

\usepackage{fancyhdr}
\usepackage{abstract}

\newtheorem{th1}{{\bf Theorem}}[section]

\newtheorem{pr1}[th1]{{\bf Proposition}}

\theoremstyle{remark}

\theoremstyle{definition}

\newtheorem{th2}{{\bf Theorem}}[subsection]
\newtheorem{le2}[th2]{{\bf Lemma}}
\newtheorem{pr2}[th2]{{\bf Proposition}}

\theoremstyle{remark}
\newtheorem{re2}[th2]{Remark}

\theoremstyle{definition}
\newtheorem{de2}[th2]{Definition}
\newtheorem{des2}[th2]{Definitions}

\theoremstyle{remark}

\theoremstyle{definition}

\DeclareGraphicsExtensions{.jpg,.mps,.pdf,.png,.gif}

\begin{document}

\title[Gysin-$(\mathbb{Z}/2\mathbb{Z})^{d}$-functors]{Gysin-$(\mathbb{Z}/2\mathbb{Z})^{d}$-functors}
\author{Dorra BOURGUIBA and Said ZARATI}
\address{Université Tunis-ElManar, Facult\'e des Sciences de Tunis, D\'epartement de Math\'ematiques.
TN-2092 Tunis, TUNISIE.}
\email{dorra.bourguiba@fst.rnu.tn}
\email{said.zarati@fst.rnu.tn}
\dedicatory{}

\subjclass{Algebraic topology}
\date{}

\maketitle
\begin{abstract} { \it Let $d \geq 1$ be an integer and $\mathcal{K}_{d}$ be a contravariant functor from  the category of subgroups of $(\mathbb{Z}/2\mathbb{Z})^{d}$ to the category of graded and finite $\mathbb{F}_{2}$-algebras. In this paper, we generalize the conjecture of G. Carlsson \cite{C3}, concerning free actions of $(\mathbb{Z}/2\mathbb{Z})^{d}$ on finite CW-complexes,  by suggesting, that if $\mathcal{K}_{d}$ is a Gysin-$(\mathbb{Z}/2\mathbb{Z})^{d}$-functor (that is to say, the functor $\mathcal{K}_{d}$ satisfies some properties, see 2.2), then we have:
$$\big(C_{d} \big): \; \underset{i \geq 0}{\sum}dim_{\mathbb{F}_{2}} \big(\mathcal{K}_{d}(0)\big)^{i} \geq 2^{d}$$
We prove this conjecture for $1 \leq d \leq 3$ and we show that, in certain cases, we get an independent proof of the following results (for $d=3$ see \cite{C4}):
\begin{center}
 If the group $(\mathbb{Z}/2\mathbb{Z})^{d}$, $ 1 \leq d \leq 3$,  acts freely and cellularly on a finite CW-complex $X$, then ${\underset{i \geq 0}{\sum}}dim_{\mathbb{F}_{2}}H^{i}(X;\; \mathbb{F}_{2}) \geq 2^{d}$.
\end{center}
}
\end{abstract}
\section{Introduction.}
Since the work of Paul A. Smith arround $1938$ \cite{Sm} (see also \cite{MB}) known as
"Smith theory" the following problem has been posed.
\begin{center}
{\bf \it $(\mathcal{P}_{d,k})$: Which group $(\mathbb{Z}/2\mathbb{Z})^{d}$ acts freely and cellularly on a
product of $k$ spheres ?}
\end{center}
The case $k = 1$, which is easy, was proved since $1935$ (\cite{Sm}, \cite{M}
and \cite{MTW}); the result is that, if $(\mathbb{Z}/2\mathbb{Z})^{d}$
acts freely on the sphere $S^{n}$ then $d \leq 1$.
\medskip\\
The case $k=2$ has been proved by A.Heller \cite{He} in $1959$
using a combinatory method which, apparently, doesn't extend to the case
of three spheres (see \cite{DV}). The result is that, if $(\mathbb{Z}/2\mathbb{Z})^{d}$
acts freely and cellularly on the product of two spheres $S^{n_{1}} \times S^{n_{2}}$, then $d \leq 2$.
\medskip\\
Some works concerning the problem $(\mathcal{P}_{d,k})$ ( such as \cite{Co}) allow to have
a generalization. The following  statement has been conjectured by Benson-Carlson \cite{BC} :
\vspace*{0.5 cm}
\begin{center}
{\bf \it ($C_{d, S}$): The group $(\mathbb{Z}/2\mathbb{Z})^{d+1}$ doesn't act
freely and cellularly on a product of $d$ spheres, $d \geq 1$}.
\end{center}
\vspace*{0.5 cm}
The conjecture ($C_{3, S}$) was proved by G. Carlsson in $1987$ \cite{C4}.
\medskip\\
Among other works concerning the conjecture ($C_{d, S}$), we can cite \cite{AB}, \cite{C2}, \cite{Han},  \cite{OY} and \cite{R}. Carlsson's work \cite{C2} concerns the case where the spheres have the same dimension and the action of the group $(\mathbb{Z}/2\mathbb{Z})^{d}$ on homology is trivial.  The work of Adem-Browder \cite{AB} concerns the case of the group $(\mathbb{Z}/p\mathbb{Z})^{d}$, $p$ an odd prime.
\medskip\\
In the middle of $1980s$ the conjecture ($C_{d, S}$) was generalized by G. Carlsson \cite{C3} (and S. Halperin \cite{Hal} for Torus) who suggest the following "Halperin-Carlsson" Conjecture (it is also called "toral rank conjecture" in some literature).
\vspace*{0.5 cm}
\begin{center}
{\bf \it ($C_{d,X}$): Let $X$ be a finite CW-complex on which the group $(\mathbb{Z}/2\mathbb{Z})^{d}$ acts freely and cellularly . Then, $\underset{i \geq 0}{\sum}dim_{\mathbb{F}_{2}}H^{i}(X;\; \mathbb{F}_{2}) \geq 2^{d}$}
\end{center}
\vspace*{0.5 cm}
In this paper we generalize the conjecture $C_{d,X}$ in the following sense which will be more precise in paragraph 2.2. Let's call a Gysin-$(\mathbb{Z}/2\mathbb{Z})^{d}$-functor a contravariant functor $\mathcal{K}_{(\mathbb{Z}/2\mathbb{Z})^{d}}$, or $\mathcal{K}_{d}$ for simplicity, from  the category of subgroups of $(\mathbb{Z}/2\mathbb{Z})^{d}$ to the category of graded, finite and  unitary $\mathbb{F}_{2}$-algebras such that:
\begin{itemize}
\item For every subgroup $W$ of $(\mathbb{Z}/2\mathbb{Z})^{d}$, the graded, finite and  unitary $\mathbb{F}_{2}$-algebra $\mathcal{K}_{d}(W)$ is non trivial and is an $H^{*}(W;\; \mathbb{F}_{2})$-algebra,
  \medskip
  \item For every subgroup $W$ of $(\mathbb{Z}/2\mathbb{Z})^{d}$ and for every $U$ a subgroup of $W$ of codimension one, there exist an exact sequence of $H^{*}(W; \;\mathbb{F}_{2})$-modules of the form:
      $$\xymatrix{...\ar[r] &\mathcal{K}_{d}(W)^{*-1} \ar[r]^{t.} & \mathcal{K}_{d}(W)^{*} \ar[r]^{\mathcal{K}_{d}(i)}
& \mathcal{K}_{d}(U)^{*} \ar[r]^{\psi} & \mathcal{K}_{d}(W)^{*} \ar[r]^{t.} & ...}$$
\end{itemize}
where
\begin{itemize}
  \item [-] $i: U \hookrightarrow W$ is the inclusion,
  \item [-] $\mathcal{K}_{d}(U)$ is an $H^{*}(W; \;\mathbb{F}_{2})$-algebra via $i^{*}: H^{*}(W; \;\mathbb{F}_{2}) \rightarrow H^{*}(U; \;\mathbb{F}_{2})$,
  \item [-] $\mathcal{K}_{d}(W)$ is an $H^{*}(W/U; \;\mathbb{F}_{2})$-algebra via $\pi^{*}: H^{*}(W/U; \;\mathbb{F}_{2}) \rightarrow H^{*}(W; \;\mathbb{F}_{2})$, $\pi: W \rightarrow W/U$ is the natural projection,
  \item [-] $H^{*}(W/U;\; \mathbb{F}_{2}) \simeq \mathbb{F}_{2}[t]$.
\end{itemize}
We propose the following conjecture:
\begin{center}
{\bf \it ($C_{d}$): Let $\mathcal{K}_{d}$ be a Gysin-$(\mathbb{Z}/2\mathbb{Z})^{d}$-functor, then: $\underset{i \geq 0}{\sum}dim_{\mathbb{F}_{2}} \big(\mathcal{K}_{d}(0)\big)^{i} \geq 2^{d}$}.
\end{center}
\vspace*{0.5 cm}
The conjecture $C_{d}$ implies the conjecture $C_{d, X}$ because if $X$ is a finite CW-complex on which the group $(\mathbb{Z}/2\mathbb{Z})^{d}$ acts freely and cellularly, then the functor $\mathcal{K}_{d}$ defined by $\mathcal{K}_{d}(W) = H^{*}_{W}(X; \; \mathbb{F}_{2})$ is a Gysin-$(\mathbb{Z}/2\mathbb{Z})^{d}$-functor whose $0^{th}$-term is $H^{*}(X; \; \mathbb{F}_{2})$.
 \medskip\\
The aim of this paper is to prove, in certain cases, the conjecture $C_{d}$ for $1 \leq d \leq 3$.
\medskip\\
The paper is structured as follows. In the second paragraph we fix the
notations and we give some properties of Gysin-$(\mathbb{Z}/2\mathbb{Z})^{d}$-functors. The third paragraph will concern the proof of the main result of this paper.
\bigskip\\
\textbf{Acknowledgements} The first author is member of the "Laboratoire de Recherche {\bf LATAO} code LR11ES12" of the Faculty of Sciences of Tunis at the  University of Tunis El-Manar.\\
The second author would like to thank "Cabinet Karray" at Avicenne for their hospitality during the preparation of this paper.
\section{On  Gysin-$(\mathbb{Z}/2\mathbb{Z})^{d}$-functors}
In this paragraph we fix some notations, introduce the Gysin-$(\mathbb{Z}/2\mathbb{Z})^{d}$-functors and give some of their properties.
\subsection{Notations} Let $V$ be an elementary abelian $2$-group that is, a group isomorphic to $(\mathbb{Z}/2\mathbb{Z})^{d},\;d \geq 1$; the integer $d$ is called the rank of $V$ and will be denoted by $d = rk(V)$. The mod. $2$ cohomology of $V$ will be simply denoted $H^{*}V$. Let's recall that $H^{*}V$ is a polynomial algebra over $ \mathbb{F}_{2}$ on  $d$ generators $t_{i}$, $1 \leq i \leq d$, of degree one. \\ We denote by  $\big(t \big)_{0}^{k} = \mathbb{F}_{2}[t] / <t^{k+1}>$ where $<t^{s}>$, $s \in \mathbb{N}$, is the ideal of $\mathbb{F}_{2}[t]$ of elements of degree $\geq s$. \\ Let $X$ be a CW-complex. Throughout this paper, the action of $V$ on $X$ will be considered  cellulary (see  \cite{TD}  [Chap. II, Sect. 1] for the notion of equivariant CW-complexes).
\subsection{Gysin-$V$-functors}
Let $V$ be an elementary abelian $2$-group of rank $\geq 1$. The set $\mathcal{W}$ of subgroups of $V$ is ordered by inclusion and then can be considered as a category. Let $\mathbb{K}_{f}$ be the category of graded, finite and unitary $\mathbb{F}_{2}$-algebras; we denote by $H^{*}V$-$\mathbb{K}_{f}$ the category of graded, finite and unitary $H^{*}V$-$\mathbb{F}_{2}$-algebras. An object of this category is a graded, finite and unitary $\mathbb{F}_{2}$-algebra $K$ equipped with a map of graded unitary $\mathbb{F}_{2}$-algebras  $H^{*}V \otimes K \rightarrow K$.
\begin{de2} A Gysin-$V$-functor is a contravariant functor
$$\mathcal{K}_{V}: \mathcal{W} \rightsquigarrow \mathbb{K}_{f}, \; W \mapsto \mathcal{K}_{V}(W) = K_{W}$$ such that:
\begin{itemize}
\item [(i)] For every subgroup $W$ of $V$, the algebra $K_{W}$ is a non trivial object of the category $H^{*}W$-$\mathbb{K}_{f}$.
\item [(ii)] For every subgroup $W$ of $V$ and for every subgroup $U$ of $W$ of codimension one, there exist an exact sequence of $H^{*}W$-modules of the form:
      $$G(U,W):\;\;\xymatrix{...\ar[r] & (K_{W})^{*-1} \ar[r]^{t.} & (K_{W})^{*} \ar[rr]^{\mathcal{K}_{V}(i)}
&& (K_{U})^{*} \ar[r]^{\psi} & (K_{W})^{*} \ar[r]^{t.} & ...}$$
\end{itemize}
where
\begin{itemize}
  \item [$\bullet$] $i: U \hookrightarrow W$ is the inclusion,
  \item [$\bullet$] $K_{U}$ is an $H^{*}W$-algebra via $i^{*}: H^{*}W \rightarrow H^{*}U$,
  \item [$\bullet$] $H^{*}(W/U) \simeq \mathbb{F}_{2}[t]$ and $t.: K_{W} \rightarrow K_{W}$ is the $H^{*}(W/U)$-structure of $K_{W}$ via the morphism $\pi^{*}: H^{*}(W/U) \rightarrow H^{*}W$ induced by the projection $\pi: W \rightarrow W/U$.
\end{itemize}
\end{de2}
{\bf 2.2.2. Vocabularies and notations}
\medskip\\
 2.2.2.1. The exact sequence $G(U,W)$ will be called the Gysin sequence associated to the subgroups $U$ and $W$ of  $V$ ($U \subseteq W$ of codimension one).\\
2.2.2.2. When the structures and morphisms are fixed, a Gysin-$V$-functor will be simply denoted $\mathcal{K}_{V}=\{K_{W},\; W \; \text{subgroup of} \;V\}$.
\medskip\\
{\bf Remark 2.2.3.} Let $\widetilde{H}^{*}(W/U)$ be the augmentation ideal of $H^{*}(W/U)$. Since $K_{W}$ is an $H^{*}(W/U)$-module via $\pi^{*}: H^{*}(W/U)\rightarrow H^{*}W$ where $\pi: W \rightarrow W/U$ is the projection; we denote
\begin{center}
$\begin{array}{lll}
                    \overline{K_{W}}^{W/U}& = &K_{W}/ \widetilde{H}^{*}(W/U).K_{W}
\medskip\\
                   & = & \mathbb{F}_{2} {\underset{\mathrm{H}^{\ast}(W/U)}
{\otimes}}K_{W} = Tor_{0}^{\mathrm{H}^{\ast}(W/U)}(K_{W},\; \mathbb{F}_{2})
  \end{array}$
\end{center}
\vspace*{0.5 cm}
The previous Gysin sequence $G(U,W)$ induces a short exact sequence of
$H^{*}U$-modules:
$$ \overline{G}(U,W):\;\; \xymatrix{0 \ar[r] &\overline{K_{W}}^{W/U}
\ar[rr]^-{\mathcal{K}_{V}(i)}&&K_{U}\ar[rr]^-{\psi} && \tau^{W/U}\big(K_{W} \big)\ar[r] & 0}
$$
where
\begin{center}
$\begin{array}{lll}
                    \tau^{W/U} \big(K_{W}\big) & = & \ker \big( t. : (K_{W})^{*}
\rightarrow (K_{W})^{*+1} \big)
\medskip\\
                   & = &Tor_{1}^{\mathrm{H}^{\ast}(W/U)}\big(K_{W},\; \mathbb{F}_{2}\big)
\end{array}$
\end{center}
\medskip
One can construct various examples of Gysin-$V$-functors; some of them are purely algebraic examples and the other comes from topology.
\medskip\\
{\bf 2.2.4. Examples }
\medskip\\
{\bf Example 2.2.4.1}. Let $K_{0}=\langle \iota, x_{1}, x_{2}, x_{4}, y_{4}, x_{5} \rangle$ be the graded, finite and unitary $\mathbb{F}_{2}$-algebra generated by six generators: $\iota$ of degree zero, $x_{i}$  of degree $i$, $i = 1, 2, 4, 5$ and $y_{4}$ of degree  $4$. These generators satisfy the following relations:
$$\left\{
  \begin{array}{ll}
    x_{j}^{2} = y_{4}^{2}=0,\; j=1, 2, 4, 5, & \hbox{} \\
    x_{1}x_{4} =  x_{1}y_{4} = 0, & \hbox{} \\
    x_{2}x_{4} =  x_{2}x_{5} = 0, & \hbox{}\\
    x_{2}y_{4} = x_{1}x_{5}. & \hbox{}
  \end{array}
\right.$$
In $K_{0}$ the elements $x_{1}x_{2}$ and $x_{1}x_{5}$ are non trivial. As an $\mathbb{F}_{2}$-vector space \\ $K_{0}=\langle \iota, x_{1}, x_{2}, x_{1}x_{2}, x_{4}, y_{4}, x_{5}, x_{1}x_{5} \rangle$.
\medskip\\
Let $H^{*}(\mathbb{Z}/2\mathbb{Z}) \cong \mathbb{F}_{2}[t]$ . We consider the graded, finite and unitary $\mathbb{F}_{2}[t]$-$\mathbb{F}_{2}$-algebra\\ $K_{\mathbb{Z}/2\mathbb{Z}}= \langle \mu, t, z_{1}, z_{2} \rangle$ generated by four generators: $\mu$ of degree zero, $t$ and $z_{1}$ of degree one and $z_{2}$ of degree two. In $K_{\mathbb{Z}/2\mathbb{Z}}$ we have the relations:
$$\left\{
  \begin{array}{ll}
    z_{1}^{2}=z_{2}^{2}=0, & \hbox{} \\
    t^{5}\mu = t^{4}z_{1} = t^{4}z_{2} = 0. & \hbox{}
  \end{array}
\right.$$
In $K_{\mathbb{Z}/2\mathbb{Z}}$ the elements $z_{1}z_{2}$, $t^{4}\mu$, $t^{3}z_{1}$,  $t^{3}z_{2}$ and $t^{3}z_{1}z_{2}$ are non trivial. We then have:
$$\left\{
                                 \begin{array}{ll}
\overline{K_{\mathbb{Z}/2\mathbb{Z}}}^{\mathbb{Z}/2\mathbb{Z}} = \langle \mu, z_{1}, z_{2}, z_{1}z_{2} \rangle \; \text{as an $\mathbb{F}_{2}$-vector space}, & \hbox{} \\
                                   \tau^{\mathbb{Z}/2\mathbb{Z}}(K_{\mathbb{Z}/2\mathbb{Z}})=\{t^{4}\mu, t^{3}z_{1}, t^{3}z_{2}, t^{3}z_{1}z_{2}\} \; \text{as an $\mathbb{F}_{2}$-vector space}. & \hbox{}
                                 \end{array}
                               \right.$$
Consider the following sequence of $\mathbb{F}_{2}$-vector spaces\\
$\xymatrix{0 \ar[r] &\overline{K_{\mathbb{Z}/2\mathbb{Z}}}^{\mathbb{Z}/2\mathbb{Z}}= \langle \mu, z_{1}, z_{2}, z_{1}z_{2} \rangle
\ar[r]^-{\sigma}&K_{0}=\langle \iota, x_{1}, x_{2}, x_{1}x_{2}, x_{4}, y_{4}, x_{5}, x_{1}x_{5} \rangle \ar[r]^-{\psi} &}$\\
\hspace*{8cm} $\xymatrix{\tau^{\mathbb{Z}/2\mathbb{Z}}(K_{\mathbb{Z}/2\mathbb{Z}})= \{t^{4}\mu, t^{3}z_{1}, t^{3}z_{2}, t^{3}z_{1}z_{2}\} \ar[r] & 0}  $\\
where
\begin{itemize}
  \item $\sigma(\mu) = \iota$, $\sigma(z_{1})=x_{1}$, $\sigma(z_{2})=x_{2}$ and $\sigma(z_{1}z_{2})=x_{1}x_{2}$.
  \item $\psi(\iota)=\psi(x_{1})=\psi(x_{2})=\psi(x_{1}x_{2})=0$
\item $\psi(x_{4})=t^{4}\mu$, $\psi(y_{4})=t^{3}z_{1}$, $\psi(x_{5})=t^{3}z_{2}$ and  $\psi(x_{1}x_{5})=t^{3}z_{1}z_{2}$.
\end{itemize}
\medskip
We verify that this sequence is exact and, by definition, that $\mathcal{K}_{\mathbb{Z}/2\mathbb{Z}}=\{ K_{0}, K_{\mathbb{Z}/2\mathbb{Z}} \}$ is a Gysin-$\mathbb{Z}/2\mathbb{Z}$-functor with $\mathcal{K}_{\mathbb{Z}/2\mathbb{Z}}(i) = \sigma$, $i: \{0\} \hookrightarrow \mathbb{Z}/2\mathbb{Z}$ is the inclusion.
\medskip\\
{\bf Example 2.2.4.2}. Let $V$ be an elementary abelian $2$-group and let $X$ be a finite  CW-complex on which the group $V$  acts freely. For every subgroup $W$ of $V$, we denote by $X_{hW}= EW \times_{W}X$ the Borel construction which is the quotient of $EW \times X$ by the diagonal action of $W$. Here $EW$ is a contractible space on which $W$ acts freely; $BW = EW/W$ is a classifying space of $W$.
The mod.$2$ cohomology of the space $X_{hW},\; H^{*}(X_{hW}) = H^{*}_{W}X $, is called the mod.$2$ equivariant cohomology of $X$. We denote by $\pi_{W}: X_{hW} \rightarrow BW$ the map induced by $X \rightarrow \{* \}$. It is clear that $H^{*}_{W}X$ is a graded $H^{*}W$-module (resp. $H^{*}V$-module) via $\pi_{W}^{*}: H^{*}W \rightarrow H^{*}_{W}X$ (resp. via $i^{*}: H^{*}V \rightarrow H^{*}W,\; \text{where}\; i: W \hookrightarrow V$ is the natural inclusion).
We verify that $H^{*}_{W}X$ is an object of the category $H^{*}W$-$\mathbb{K}_{f}$ that is a graded, finite and unitary $\mathbb{F}_{2}$-algebra equipped with a map of graded unitary $\mathbb{F}_{2}$-algebras  $H^{*}W \otimes H^{*}_{W}X \rightarrow H^{*}_{W}X$.
\smallskip\\
The contravariant functor $\mathcal{K}_{V}: \mathcal{W} \rightsquigarrow \mathbb{K}_{f},\; W \mapsto \mathcal{K}_{V}(W) = H^{*}_{W}X$ is a Gysin-$V$-functor because:
\begin{itemize}
\item [(i)] For every subgroup $W$ of $V$, $H^{*}_{W}X$ is a non trivial object of the category $H^{*}W$-$\mathbb{K}_{f}$.
\item [(ii)] Let $W \subseteq V$ be a subgroup  and let $U \subset W$ be a subgroup of codimension $1$. The inclusion $i: U \hookrightarrow W$ induces the following two sheets covering: $ W/U \cong \mathbb{Z}/2\mathbb{Z} \rightarrow X_{hU} \overset{i} \rightarrow X_{hW}$ with
$B(\pi) \circ \pi_{W}:  X_{hW} \rightarrow BW \rightarrow B(W/U)$ as a classifying map, $\pi: W \rightarrow W/U$ is the natural projection.\\
Let $H^{*}(W/U) \cong \mathbb{F}_{2}[t]$, we also denote by $t$ the non
trivial element $(B(\pi) \circ \pi_{W})^{*}(t)$ of $  H^{1}_{W}X$. The Gysin exact sequence
associated to the previous covering is the following exact sequence of $H^{*}W$-modules:
\begin{center}
$\xymatrix{...\ar[r]&H_{W}^{*-1}X\ar[r]^{t.}&H_{W}^{*}X\ar[r]^{i^{*}}
&H^{*}_{U}X\ar[r]^{tr}&H_{W}^{*}X\ar[r]^{t.}&...}$
\end{center}
where $tr$ is the the transfer (\cite{Sp}, \cite{Z}) and for $ x \in H_{W}^{\ast} X$, $t.x =
(B(\pi) \circ \pi_{W})^{\ast}(t) \smile x$.
\end{itemize}
\medskip
This shows that $\mathcal{K}_{V}=\{K_{W} = H_{W}^{\ast}X,\; W \; \text{a subgroup of} \;V\}$ is a Gysin-$V$-functor. This example comes from "topology" via the equivariant cohomology of a free action of $V$  on a finite CW-complex $X$.
\subsection{Some properties of Gysin-$V$-functors}\leavevmode\par \noindent
Let's recall some definitions and fix some notations. Let $E$ be a finite graded $\mathbb{F}_{2}$-vector space.
\medskip\\
$\bullet$ We denote by $\parallel E \parallel$ the {\it norm} of $E$ which is the maximum of the set $ \{k \in \mathbb{N}, \; E^{k} \neq 0\}$.
\medskip\\
$\bullet$ Let $V$ be an elementary abelian $2$-group. If $E$ is an $H^{*}V$-module and $x \in E$, we denote by $\langle x \rangle_{V}$ the sub-$H^{*}V$-module of $E$ generated by the element $x$.
\begin{des2} (i) A finite graded $\mathbb{F}_{2}$-vector space $E$ is called:
\medskip\\
\hspace*{1cm} (i-a)  {\it connected} if $E^{0} \cong \mathbb{Z}/2\mathbb{Z}$.\\
\hspace*{1cm} (i-b) {\it bi-connected} if:
$\left\{
  \begin{array}{ll}
    \text{$E$ is connected}:\; E^{0} \cong \mathbb{Z}/2\mathbb{Z},  & \hbox{} \\
    \text{and}\medskip \\
 E^{\parallel E \parallel} \cong \mathbb{Z}/2\mathbb{Z}. & \hbox{}
  \end{array}
\right.$
\medskip\\
(ii) A Gysin-$V$-functor $\mathcal{K}_{V} = \{K_{W},\; W \; \text{subgroup of} \;V\}$ is {\it connected} (resp. {\it bi-connected}) if $K_{0}$ is {\it connected} (resp. $K_{0}$ is {\it bi-connected}).
\medskip\\
(iii) A finite CW-complex $X$ is {\it bi-connected} if $H^{\ast}X $ is {\it bi-connected}.
\end{des2}
We have the following property of Gysin-$V$-functors.
\begin{le2} Let $V$ be an elementary abelian $2$-group and let $\mathcal{K}_
{V}=\{K_{W},\; W \; \text{subgroup of} \;V\}$ be a bi-connected Gysin-$V$-functor. Then, for each subgroup $W$ of $V$, the graded finite $\mathbb{F}_{2}$-algebra $K_{W}$ is bi-connected and we have $\parallel K_{W} \parallel = \parallel K_{0} \parallel$.
\end{le2}
\begin{proof} The proof is by induction on the rank of the subgroup of $V$. Let $\mathcal{K}_{V}$ be a bi-connected Gysin-$V$-functor and let $U \subseteq V$ be a subgroup of rank one. The Gysin exact sequence of graded $\mathbb{F}_{2}$-vector spaces:
$$ \overline{G}(0, U):\;\; \xymatrix{0 \ar[r] &\overline{K_{U}}^{U}
\ar[rr]^-{\mathcal{K}_{V}(i)}&&K_{0} \ar[rr]^-{\psi} && \tau^{U}(K_{U})\ar[r] & 0}
$$
shows that:
\medskip\\
{\bf 2.3.2.1.} In degree zero, we have: $(\mathcal{P}_{0}): \;  \mathbb{Z}/2\mathbb{Z} \cong \big(K_{0}\big)^{0} \cong \big(\overline{K_{U}}^{U}\big)^{0} \oplus \big(\tau^{U}(K_{U})\big)^{0}$.
\medskip\\
Since $\big(\tau^{U}(K_{U}) \big)^{0} \subseteq \big(K_{U}\big)^{0} = \big( \overline{K_{U}}^{U} \big)^{0}$, then $\big( \overline{K_{U}}^{U} \big)^{0} = 0$ implies $\big(\tau^{U}(K_{U}) \big)^{0} = 0$. This contradicts the equality $(\mathcal{P}_{0})$. Then we deduce that: $\big(\tau^{U}(K_{U}) \big)^{0} = 0$ and  $\mathbb{Z}/2\mathbb{Z} \cong \big( \overline{K_{U}}^{U} \big)^{0} \cong \big(K_{U}\big)^{0}$. This shows that $K_{U}$ is connected.
\medskip\\
{\bf 2.3.2.2.} In degree $\parallel K_{0} \parallel$, we have: $(\mathcal{P}_{\parallel K_{0} \parallel}): \; \mathbb{Z}/2\mathbb{Z} \cong \big(K_{0}\big)^{\parallel K_{0} \parallel} \cong \big(\overline{K_{U}}^{U}\big)^{\parallel K_{0} \parallel} \oplus \big(\tau^{U}(K_{U})\big)^{\parallel K_{0} \parallel}$. Since $K_{U}$ is a graded finite $H^{*}U$-module, we have: $\parallel K_{U} \parallel = \parallel \tau^{U}(K_{U}) \parallel$. The following inequalities follow: $ \parallel \overline{K_{U}}^{U} \parallel  \leq \parallel K_{U} \parallel = \parallel \tau^{U}(K_{U}) \parallel \leq \parallel K_{0} \parallel$. We deduce from $(\mathcal{P}_{\parallel K_{0} \parallel})$ that: $ \big( \overline{K_{U}}^{U} \big)^{\parallel K_{0} \parallel} = 0$ and $\big(\tau^{U}(K_{U})\big)^{\parallel K_{0} \parallel} \cong  \mathbb{Z}/2\mathbb{Z}$. This shows that: \\ $\parallel K_{U} \parallel = \parallel \tau^{U}(K_{U}) \parallel \geq \parallel K_{0} \parallel$. So, we have the equality: $\parallel K_{U} \parallel = \parallel \tau^{U}(K_{U}) \parallel = \parallel K_{0} \parallel$.
\medskip \\
We proved that $K_{U}$ is bi-connected and  $ \parallel K_{U} \parallel = \parallel K_{0} \parallel$.
\medskip \\
The lemma holds by induction on the rank of subgroups of $V$ using the same method.
\end{proof}
Let $E$ be a graded finite $\mathbb{F}_{2}$-vector space. We denote by $d(E) = \underset{i \geq 0}{\sum}dim_{\mathbb{F}_{2}}E^{i}$ the (total) dimension of $E$. We have:
\begin{pr2} Let $V$ be an elementary abelian $2$-group and $\mathcal{K}_{V}=\{K_{W},\; W \; \text{subgroup of} \;V\}$ be a Gysin-$V$-functor. Then, the dimension of $K_{0}$ is even:
$d(K_{0}) \equiv 0 \; (mod.\; 2)$.
\end{pr2}
\begin{proof} Let $U \subset V$ be a subgroup of rank one, then the Gysin exact sequence
$$ \overline{G}(0, U):\;\; \xymatrix{0 \ar[r] &\overline{K_{U}}^{U}
\ar[rr]^-{\mathcal{K}_{V}(i)}&&K_{0} \ar[rr]^-{\psi} && \tau^{U}(K_{U})\ar[r] & 0}
$$
Shows that $d(K_{0}) = d\big(\overline{K_{U}}^{U}\big) + d\big(\tau^{U}(K_{U})\big)$. The propostion 2.3.3 is a consequence of the following lemma.
\end{proof}
\begin{le2} Let $U$ be an elementary abelian $2$-group of rank one and let $M$ be a graded finite $H^{*}U$-module. We have: $d\big(\overline{M}^{U}\big) = d\big(\tau^{U}(M)\big)$.
\end{le2}
\begin{proof} The proof is by induction on the dimension of the finite $\mathbb{F}_{2}$-vector space $\overline{M}^{U}$.
\end{proof}
\begin{re2} Here is an example of application of the previous lemma. Let $U_{i}$, $i=1, 2$, be an elementary abelian $2$-group of rank one, let $V = U_{1} \oplus U_{2}$ and let $H^{*}U_{i} \cong \mathbb{F}_{2}[t_{i}], \; i=1,2$. Let $M$ be a graded finite $H^{*}V$-module, $x_{1}$ and $x_{2}$ two elements of $M$ such that the finite $ \mathbb{F}_{2}[t_{1}]$-modules
\smallskip\\
 $ \overline{M}^{U_{2}}$ and $\tau^{U_{2}}(M)$ are monogenic:
$\left\{
  \begin{array}{ll}
    \overline{M}^{U_{2}}  \cong \big(t_{1} \big)_{0}^{k_{1}}x_{1}, & \hbox{}
\medskip\\
   \tau^{U_{2}}(M) \cong  \big(t_{1} \big)_{0}^{k_{2}}x_{2}. & \hbox{}
  \end{array}
\right.$ (see notations 2.1).
\medskip\\
In this case we have $d \big( \overline{M}^{U_{2}} \big) = k_{1} + 1$ and $d \big( \tau^{U_{2}}(M) \big) = k_{2} + 1$. The lemma 2.3.4, applied for the $H^{*}U_{2}$-module $M$, implies the equality: $k_{1} = k_{2}$.
\end{re2}
\subsection{On the extension of Gysin-$V$-functors}
Let $V$ be an elementary abelian $2$-group, $V'$ a subgroup of $V$ and $\mathcal{K}_{V}=  \{K_{W},\; W \; \text{subgroup of} \;V \}$ be a Gysin-$V$-functor. Then, $\mathcal{K}_{V'}=\{K_{W},\; W \; \text{subgroup of} \;V' \}$ is a Gysin-$V'$-functor called a "sub-Gysin-functor" of $\mathcal{K}_{V}$. We say also that $\mathcal{K}_{V}$ is an extension of the Gysin-$V'$-functor $\mathcal{K}_{V'}$.
\medskip\\
It is interesting to know when a Gysin-$V$-functor extends because this question is related to the extension of a free action of the group $V$ on a finite CW-complex. We have.
\begin{pr2} Let $V$ be an elementary abelian $2$-group and let $\mathcal{K}_{V}=\{K_{W},\; W \; \text{subgroup of} \;V\}$ be a Gysin-$V$-functor such that $\overline{K_{V}}^{V} \cong \mathbb{Z}/2\mathbb{Z}$, then $\mathcal{K}_{V}$ doesn't extend.
\end{pr2}
\begin{proof} Suppose that the Gysin-$V$-functor $\mathcal{K}_{V}$ extends to $\mathcal{K}_{H}$, $H = V \oplus \mathbb{Z}/2\mathbb{Z}$ then, the Gysin exact sequence of graded finite $H^{*}V$-modules:
$$ \overline{G}(V, H):\;\; \xymatrix{0 \ar[r] &\overline{K_{H}}^{H/V}
\ar[rr]^-{\mathcal{K}_{H}(i)}&&K_{V} \ar[rr]^-{\psi} && \tau^{H/V}(K_{H})\ar[r] & 0}
$$
Shows that $$\begin{array}{lll}
                   \big (K_{V} \big)^{0} & \cong & \big ( \overline{K_{H}}^{H/V} \big)^{0} \oplus \big ( \tau^{H/V}(K_{H}) \big)^{0}
\medskip\\
                   & \cong & \big (K_{H} \big)^{0} \oplus \big ( \tau^{H/V}(K_{H}) \big)^{0}
\end{array}$$
Since $\big ( \tau^{H/V}(K_{H}) \big)^{0} \subseteq \big (K_{H} \big)^{0}$ and $ \big (K_{V} \big)^{0} \cong \mathbb{Z}/2\mathbb{Z}$ because
$\overline{K_{V}}^{V} \cong \mathbb{Z}/2\mathbb{Z}$, we deduce that the morphism $ \xymatrix{\big ( \overline{\overline{K_{H}}^{H/V} \big)}^{V}
\ar[rr]^-{\mathcal{K}_{H}(i)}&&\overline{K_{V}}^{V} \cong \mathbb{Z}/2\mathbb{Z} }$ is an epimorphism and $\overline{\tau^{H/V}(K_{H})}^{V} = 0$. This leads to a contradiction because the $H^{*}(H)$-module $K_{H}$ is non trivial and finite.
\end{proof}
\begin{pr2} Let $V$ be an elementary abelian $2$-group, $\mathcal{K}_{V}$ a bi-connected Gysin-$V$-functor and  $\iota_{V}$ the unit of the $\mathbb{F}_{2}$-algebra $K_{V}$: $(K_{V})^{0} \cong \mathbb{Z}/2\mathbb{Z} = < \iota_{V} >$.\\ If the norm of $K_{V}$ is equal to the norm of its sub-$H^{*}V$-module generated by $\iota_{V}$:\\ $\parallel K_{V} \parallel = \parallel  \langle \iota_{V} \rangle_{V} \parallel$, then $\mathcal{K}_{V}$ doesn't extend.
\end{pr2}
\begin{proof} Let $\mathcal{K}_{V}=\{K_{W},\; W \; \text{subgroup of} \;V\}$ be a  bi-connected Gysin-$V$-functor and suppose that $\mathcal{K}_{V}$ extends to $\mathcal{K}_{H}$, where $H = V \oplus \mathbb{Z}/2\mathbb{Z}$. By the lemma 2.3.2, the Gysin-$H$-functor $\mathcal{K}_{H}$ is bi-connected and we have: $ \parallel  K_{V} \parallel = \parallel K_{H} \parallel$ and $(K_{H})^{0} \cong \mathbb{Z}/2\mathbb{Z} = < \iota_{H} >$. Since the map $\mathcal{K}_{H}(i): K_{H} \rightarrow K_{V}$, induced by the inclusion of $V$ in $H$, is a map of unitary-(connected)-$\mathbb{F}_{2}$-algebras then,  $\mathcal{K}_{H}(i)(\iota_{H})=\iota_{V}$.
\medskip\\
Let's denote by $j: \langle \iota_{H} \rangle_{H} \hookrightarrow K_{H}$ the natural inclusion. We have the following commutative diagram whose second line is the short Gysin exact sequence $\overline{G}(V, H)$ of $H^{*}V$-modules.
$$\xymatrix{ &Im(\overline{j}^{H/V}) \ar@{^{(}->}[d]\ar@{->>}[rr]^-{\mathcal{K}_{H}(i)_{\mid}}&&\langle \iota_{V} \rangle_{V} \ar@{^{(}->}[d]& \\
0 \ar[r] &\overline{K_{H}}^{H/V} \ar[rr]^-{\mathcal{K}_{H}(i)}&&K_{V} \ar[rr]^-{\psi} && \tau^{H/V}(K_{H})\ar[r] & 0}
$$
This shows that the sub-$H^{*}V$-module $Im(\overline{j}^{H/V})$ of $\overline{K_{H}}^{H/V}$ is isomorphic to sub-$H^{*}V$-module $ \langle \iota_{V} \rangle_{V}$ of $K_{V}$ \big($\mathcal{K}_{H}(i)_{\mid}: Im(\overline{j}^{H/V}) \rightarrow \langle \iota_{V} \rangle_{V}$ is an isomorphism\big). This implies the inequality between norms:  $\parallel \overline{K_{H}}^{H/V} \parallel \geq  \parallel Im(\overline{j}^{H/V}) \parallel = \parallel \langle \iota_{V} \rangle_{V} \parallel = \parallel K_{V} \parallel$.
\medskip\\
Since $\overline{K_{H}}^{H/V}$ is a sub-$H^{*}V$-module of $K_{V}$, we have that $\parallel \overline{K_{H}}^{H/V} \parallel \leq \parallel K_{V} \parallel$. So we have the equality: $\parallel \overline{K_{H}}^{H/V} \parallel = \parallel K_{V} \parallel$.
\medskip\\
The Gysin exact sequence $\overline{G}(V, H)$ of $H^{*}V$-modules:
$$\xymatrix{ 0 \ar[r] &\overline{K_{H}}^{H/V} \ar[rr]^-{\mathcal{K}_{H}(i)}&&K_{V} \ar[rr]^-{\psi} && \tau^{H/V}(K_{H})\ar[r] & 0}$$
implies the following isomorphism: $\big ( K_{V} \big )^{\parallel K_{V} \parallel} \cong \big ( \overline{K_{H}}^{H/V} \big )^{\parallel K_{V} \parallel} \oplus \big ( \tau^{H/V}(K_{H}) \big )^{\parallel K_{V} \parallel}$.
\medskip\\
Since the Gysin-$V$-functor  $\mathcal{K}_{V}$ is bi-connected, $\big ( K_{V} \big )^{\parallel K_{V} \parallel} \cong \mathbb{Z}/2\mathbb{Z}$, and $\parallel \overline{K_{H}}^{H/V} \parallel = \parallel K_{V} \parallel$, we deduce from the previous isomorphism:
\medskip
$$\left\{
  \begin{array}{ll}
   (i) \; \big ( \overline{K_{H}}^{H/V} \big )^{\parallel K_{V} \parallel} \cong \big ( K_{V} \big )^{\parallel K_{V} \parallel} \cong \mathbb{Z}/2\mathbb{Z} , & \hbox{}
\medskip\\
  (ii) \; \big ( \tau^{H/V}(K_{H}) \big )^{\parallel K_{V} \parallel} = 0. & \hbox{}
  \end{array}
\right.$$
By proposition 2.3.2, $\parallel K_{V} \parallel = \parallel K_{H} \parallel$; since $K_{H}$ is a graded, finite and non trivial $H^{*}(H/V)$-module, then $\big ( \tau^{H/V}(K_{H}) \big )^{\parallel K_{V} \parallel} \neq 0$. This contradicts the point $(ii)$.
\end{proof}
\section{The main result}
 Let $V$ be an elementary abelian $2$-group of rank $d$ and $\mathcal{K}_{V}=\{K_{W},\; W \; \text{subgroup of} \;V\}$ be a Gysin-$V$-functor. Let's denote $d(K_{0})=\Sigma_{i \geq 0}dim_{\mathbb{F}_{2}} \big( K_{0} \big)^{i}$  the total dimension of the graded finite $\mathbb{F}_{2}$-vector space $K_{0}$. \\
The main result of this paper is to show, in certain cases, that  $d(K_{0})$ is related to the rank of the group $V$, as suggested by the conjecture $(C_{d})$, $d(K_{0}) \geq 2^{rk(V)}$. \\
More precisely, we have:
\begin{th1} Let $V$ be an elementary abelian $2$-group and $\mathcal{K}_{V}=\{K_{W},\; W \; \text{subgroup of} \;V\}$ be a Gysin-$V$-functor. Then,\\
(i) For $rk(V)=1$, $d(K_{0}) \geq 2$ so the conjecture $(C_{1})$ holds.\\
(ii) For $rk(V)=2$, if the Gysin-$V$-functor $\mathcal{K}_{V}$ is connected, we have the inequality: $d(K_{0}) \geq 4$, so the conjecture $(C_{2})$ holds.\\
(ii) For $rk(V)=3$, if the Gysin-$V$-functor $\mathcal{K}_{V}$ is  bi-connected, we have the inequality:\\ $d(K_{0}) \geq 8$, so the conjecture $(C_{3})$ holds.
\end{th1}
As an application of this theorem we get an independent proof of the results concerning $(C_{d,X})$ for $d \leq 3$.
\begin{pr1} Let $V$ be an elementary abelian $2$-group and let $X$ be a finite CW-complex on which the group $V$ acts freely. Then, \\
(i) For $rk(V)=1$, we have: $d(H^{*}X) \geq 2$.\\
(ii) For $rk(V)=2$ and $X$ connected, we have: $d(H^{*}X) \geq 4$.\\
(ii) For $rk(V)=3$ and $X$ bi-connected, we have: $d(H^{*}X) \geq 8$.
\end{pr1}
\begin{proof} Let $V$ be an elementary abelian $2$-group and let $X$ be a finite CW-complex on which the group $V$ acts freely. By the example 2.2.4.2, the contravariant functor $\mathcal{K}_{V}: \mathcal{W} \rightsquigarrow \mathbb{K}_{f},
\; W \mapsto \mathcal{K}_{V}(W) = H^{*}_{W}X$ is a Gysin-$V$-functor whose $0^{th}$-term $K_{0}=
\mathcal{K}_{V}(0) = H^{*}X$.
\end{proof}
Let $S^{n}$, $n \geq 1$, be the standard unit sphere in $\mathbb{R}^{n+1}$, then the product $S^{n_{1}} \times...\times S^{n_{k}}$, $k \geq 1$, is a bi-connected CW-complex. By the proposition 3.2, if an elementary abelian $2$-group $V$,  $1 \leq rk(V) \leq 3$, acts freely on a product of $k$  spheres then, $k \geq rk(V)$.
\subsection{proof of theorem 3.1}\leavevmode\par \noindent
To prove theorem $3.1$ we consider the following three cases:
\medskip\\
{\bf 3.1.1. { \it The case $rk(V)=1$.}}
\medskip\\
The proposition 2.3.3 shows that if $\mathcal{K}_{V}=\{K_{W},\; W \; \text{subgroup of} \;V\}$ is a Gysin-$V$-functor, $rk(V)=1$, then $d(K_{0}) \equiv 0\; (mod.\; 2)$. This implies that $d(K_{0}) \geq 2$ because the graded $\mathbb{F}_{2}$-vector space $K_{0}$ is not trivial.
\medskip\\
{\bf 3.1.2. { \it The case $rk(V)=2$.}}
\medskip\\
Let $\mathcal{K}_{V}=\{K_{W},\; W \; \text{subgroup of} \;V\}$ be a Gysin-$V$-functor, $rk(V)=2$,  and suppose that $d(K_{0}) < 4$. Since $d(K_{0}) \equiv 0\; (mod.\; 2)$ (see proposition 2.3.3 ) and $K_{0}$ non trivial, we deduce that $d(K_{0}) = 2$.
\medskip\\
Let $U \subseteq V$ be a subgroup of rank one and consider the short exact sequence of graded $\mathbb{F}_{2}$-vector spaces associated to the couple $(\{0 \} \subset U)$ of subgroups of $V$
$$ \overline{G}(\{0\}, U):\;\; \xymatrix{0 \ar[r] &\overline{K_{U}}^{U}
\ar[rr]^-{\mathcal{K}_{V}(i)}&&K_{0} \ar[rr]^-{\psi} && \tau^{U}(K_{U})\ar[r] & 0}
$$
$i: \{0\} \hookrightarrow U$ denotes the inclusion. This shows that: $d(K_{0}) = 2 = d(\overline{K_{U}}^{U}) + d(\tau^{U}(K_{U}))$. The lemma 2.3.4 implies that: $d(\overline{K_{U}}^{U}) =1$. Since the Gysin-$V$-functor $\mathcal{K}_{V}$ is {\bf connected}, we have: $\mathbb{Z}/2\mathbb{Z} \cong (K_{U})^{0} \cong (\overline{K_{U}}^{U})^{0} \cong \overline{K_{U}}^{U}$.\\
The proposition 2.4.1 shows that, in this case, the Gysin-$U$-functor $\mathcal{K}_{U}=\{K_{W},\; W \; \text{subgroup of} \;U\}$ can not extend to $\mathcal{K}_{V}$. This leads to a contradiction.
\medskip\\
{\bf 3.1.3. { \it The case $rk(V)=3$.}}
\medskip\\
Let $\mathcal{K}_{V}=\{K_{W},\; W \; \text{subgroup of} \;V\}$ be a bi-connected Gysin-$V$-functor, $rk(V)= 3$, and suppose that $d(K_{0}) < 8$. Since $d(K_{0}) \equiv 0\; (mod.\; 2)$ (see proposition 2.3.3 ) and $K_{0}$ non trivial, then we have three possibility: $d(K_{0}) = 2$, $d(K_{0}) = 4$ and $d(K_{0}) = 6$. We will show that the three cases $d(K_{0}) = 2$, $d(K_{0}) = 4$ and $d(K_{0}) = 6$ are impossible. Let $U_{i}, 1 \leq i \leq 3$, be a rank one subgroup of $V$ such that: $V \cong U_{1} \oplus U_{2} \oplus U_{3}$.
\medskip\\
{\bf 3.1.3.1} The case $d(K_{0}) = 2$ is impossible by the previous case 3.1.2. We proved in 3.1.2 that if $d(K_{0}) = 2$, then $K_{0}$ can't be the $0^{th}$-term of a Gysin-$E$-functor with $E$ an elementary abelian $2$-group of rank $2$ and a fortiori of rank $\geq 2$.
\medskip\\
{\bf 3.1.3.2} Suppose that $d(K_{0}) = 4$. The Gysin exact sequence of graded finite  $\mathbb{F}_{2}$-vector spaces
$$ \overline{G}(\{0\}, U_{1}):\;\; \xymatrix{0 \ar[r] &\overline{K_{U_{1}}}^{U_{1}}
\ar[rr]^-{\mathcal{K}_{V}(i_{1})}&&K_{0} \ar[rr]^-{\psi} && \tau^{U_{1}}(K_{U_{1}})\ar[r] & 0}
$$
($i_{1}: \{0\} \hookrightarrow U_{1}$ is the inclusion), shows that $d(K_{0}) = 4 = d\big(\overline{K_{U_{1}}}^{U_{1}} \big) + d\big(\tau^{U_{1}}(K_{U_{1}}) \big)$. The lemma 2.3.4 implies that $d(\overline{K_{U_{1}}}^{U_{1}}) = 2$, that is: $\overline{K_{U_{1}}}^{U_{1}} \cong \langle \overline{\iota_{1}}, \overline{g_{1}} \rangle$ is the $\mathbb{F}_{2}$-vector space generated by two generators $\overline{\iota_{1}}$ and $\overline{g_{1}}$ where $\iota_{1} \in \big(K_{U_{1}}\big)^{0} \cong \mathbb{Z}/2\mathbb{Z}$ is the unit and $g_{1} \in \big(K_{U_{1}}\big)^{k_{1}}$, $k_{1} \geq 1$.
\medskip\\
Since $\mathcal{K}_{U_{1}}$ is a sub-Gysin-functor of $\mathcal{K}_{U_{1} \oplus U_{2} }$ whose $0^{th}$-term $K_{0}$ is bi-connected, then by 2.4.2, the norm of $K_{U_{1}}$ is bigger than the norm of the sub-$H^{*}U_{1}$-module generated by $\iota_{1}$. We have:
$$\begin{array}{lll}
                  \parallel K_{U_{1}} \parallel  & = & \parallel \langle g_{1} \rangle_{U_{1}}  \parallel
\medskip\\
                 & > & \parallel \langle \iota_{1} \rangle_{U_{1}}  \parallel
\end{array}$$
This shows, in particular, that we have an isomorphism of $H^{*}U_{1}$-modules:
$$K_{U_{1}} \cong \langle \iota_{1} \rangle_{U_{1}} \oplus \langle  g_{1} \rangle_{U_{1}}$$
The Gysin exact sequence $\overline{G}(U_{1}, U_{1} \oplus U_{2})$ of $H^{*}U_{1}$-modules
$$ \xymatrix{0 \ar[r] &\overline{K_{U_{1} \oplus U_{2}}}^{U_{2}}
\ar[rr]^-{\mathcal{K}_{U_{1} \oplus U_{2}}(j_{1})}&&K_{U_{1}}\cong \langle \iota_{1} \rangle_{U_{1}} \oplus \langle  g_{1} \rangle_{U_{1}} \ar[rr]^-{\psi} && \tau^{U_{2}}(K_{U_{1} \oplus U_{2}})\ar[r] & 0}
$$
$j_{1}: U_{1} \hookrightarrow U_{1} \oplus U_{2}$ denotes the natural inclusion, shows that:
$$\left\{
              \begin{array}{ll}
               (i) \; \overline{K_{U_{1} \oplus U_{2}}}^{U_{2}} \cong \langle \iota_{1}  \rangle_{U_{1}} ,
\medskip\\
               (ii) \; \tau^{U_{2}}(K_{U_{1} \oplus U_{2}}) \cong \langle \psi(g_{1}) \rangle_{U_{1}}.
              \end{array}
            \right.$$
\medskip
The point $(i)$ implies that $\overline{K_{U_{1} \oplus U_{2}}}^{U_{1} \oplus U_{2}} \cong \mathbb{Z}/2\mathbb{Z}$ and the proposition 2.4.1 shows the contradiction since the Gysin-$V$-functor $\mathcal{K}_{V}$ extends $\mathcal{K}_{U_{1} \oplus U_{2}}$ ($V \cong U_{1} \oplus U_{2} \oplus U_{3}$).
\medskip\\
{\bf 3.1.3.3} Suppose that $d(K_{0}) = 6$. To show a contradiction, in this case, we will analyse the graded, finite and unitary $H^{*}W $-$ \mathbb{F}_{2}$-algebras $K_{W}$ for $W = U_{1}$ and $W = U_{1} \oplus U_{2}$.
\medskip\\
{\bf I1. Informations on $K_{U_{1}}$}.
\medskip\\
By the same previous method, using the Gysin exact sequence
\begin{center}
$ \overline{G}(\{0\}, U_{1}):\;\; \xymatrix{0 \ar[r] &\overline{K_{U_{1}}}^{U_{1}}
\ar[rr]^-{\mathcal{K}_{V}(i_{1})}&&K_{0} \ar[rr]^-{\psi} && \tau^{U_{1}}(K_{U_{1}})\ar[r] & 0},
$
\end{center}
we show that $d(\overline{K_{U_{1}}}^{U_{1}}) = 3$ that is: $\overline{K_{U_{1}}}^{U_{1}} \cong \langle \overline{\iota_{1}}, \overline{g_{1}}, \overline{g_{2}} \rangle$ is the $\mathbb{F}_{2}$-vector space generated by three generators $\overline{\iota_{1}}$, $\overline{g_{1}}$ and $\overline{g_{2}}$:  $\iota_{1} \in \big(K_{U_{1}}\big)^{0} \cong \mathbb{Z}/2\mathbb{Z}$ is the unit, $g_{1} \in \big(K_{U_{1}}\big)^{k_{1}}$, $k_{1} \geq 1$ and $g_{2} \in \big(K_{U_{1}}\big)^{k_{2}}$, $k_{2} \geq 1$.
\smallskip\\
Since the bi-connected $U_{1}$-Gysin functor $\mathcal{K}_{U_{1}}$ extends, then by  proposition 2.4.2, the norm of the graded finite $\mathbb{F}_{2}$-vector space $K_{U_{1}}$ is reached as the norm of a sub-$\mathbb{F}_{2}$-vector space generated by a generator different of $\iota_{1}$, for example $g_{1}$. We have: $\parallel \langle \iota_{1} \rangle_{U_{1}} \parallel < \parallel K_{U_{1}} \parallel = \parallel \langle g_{1} \rangle_{U_{1}} \parallel$. \\
We verify then that we have a short exact sequence of $H^{*}U_{1}$-modules of the form:
\begin{center}
$ \big(E(U_{1})\big): \; \xymatrix{0 \ar[r] & \langle \iota_{1} \rangle_{U_{1}} \oplus \langle g_{1} \rangle_{U_{1}}
\ar[r]&K_{U_{1}} \ar[r]& C_{U_{1}} \ar[r] & 0}$
\end{center}
where $C_{U_{1}}$ is a graded finite monogenic $H^{*}U_{1}$-module generated by the element $g_{2}$.
\smallskip\\
In refers to 2.1, let $H^{*}U_{i} \cong \mathbb{F}_{2}[t_{i}]$, $i=1, 2$, the polynomial algebra over $\mathbb{F}_{2}$ on one generator $t_{i}$ of degree one, $<t^{s}>$, $s \in \mathbb{N}$, be the ideal of $\mathbb{F}_{2}[t]$ of elements of degree $\geq s$ and\\ $ \big(t \big)_{0}^{k} = \mathbb{F}_{2}[t] / <t^{k+1}> $.
\smallskip\\
With these notations we have:
\smallskip
\begin{itemize}
  \item [\bf I1.1]  $\langle \iota_{1} \rangle_{U_{1}} \cong \big(t_{1} \big)_{0}^{n_{1}}\iota_{1}$, $n_{1} \geq 1$.
\smallskip
  \item [\bf I1.2]  $C_{U_{1}} \cong \big(t_{1} \big)_{0}^{l_{1}}g_{2}$ with \mathversion{bold} $l_{1} \leq n_{1} $ \mathversion{normal}
 because, in the graded finite unitary $\mathbb{F}_{2}$-algebra $K_{U_{1}}$, we have: $g_{2} = \iota_{1}.g_{2}$. This implies that: $t_{1}^{s}g_{2} = (t_{1}^{s}\iota_{1}).g_{2}$, $s \in \mathbb{N}$.
\end{itemize}
\medskip
{\bf I1.3 Remark}. In $I1.1$ the integer $n_{1}$ is $\geq 1$ because if not $n_{1} = 0$ which means that: $t_{1}.\iota_{1} = 0$. This implies that $\iota_{1} \in \big(\tau^{U_{1}}(K_{U_{1}})\big)^{0}$. Since $\big(K_{0}\big)^{0} \cong \big(\overline{K_{U}}^{U}\big)^{0} \oplus \big(\tau^{U}(K_{U})\big)^{0}$, we get a contradiction with $K_{0}$ connected: $ \big(K_{0}\big)^{0} \cong \mathbb{Z}/2\mathbb{Z}$.
\medskip\\
{\bf I2. Informations on $K_{U_{1} \oplus U_{2}}$}.
\medskip\\
Let $\iota_{1,2}$ be the unit of the graded $\mathbb{F}_{2}$-algebra $K_{U_{1} \oplus U_{2}}$ and consider the short exact sequence of \\ $H^{*}(U_{1} \oplus U_{2})$-modules:
\begin{center}
$ \big(E(U_{1} \oplus U_{2})\big): \; \xymatrix{0 \ar[r] &  \langle \iota_{1,2} \rangle_{U_{1} \oplus U_{2}}
\ar[r]^{j}&K_{U_{1} \oplus U_{2}} \ar[r]& C_{U_{1} \oplus U_{2}} \ar[r] & 0}$
\end{center}
We have the following commutative diagram, {\bf (D)}, of $H^{*}U_{1}$-modules:
$$   \xymatrix{&0 \ar[d]&& 0 \ar[d]&& 0 \ar[d] & \\
0 \ar[r] & Im (\overline{j}^{U_{2}})
\ar[rr] \ar[d]&& \langle \iota_{1} \rangle_{U_{1}} \oplus \langle g_{1} \rangle_{U_{1}} \ar[rr]^{\psi_{\mid}}\ar[d]&& \langle \psi(g_{1}) \rangle_{U_{1}} \ar[r]\ar[d] & 0 \\
 0 \ar[r] & \overline{K_{U_{1} \oplus U_{2}}}^{U_{2}}
\ar[rr]^{\mathcal{K}_{U_{1} \oplus U_{2}}(i_{1})}\ar[d]&&K_{U_{1}} \ar[rr]^{\psi}\ar[d]&& \tau^{U_{2}}\big(K_{U_{1} \oplus U_{2}} \big) \ar[r]\ar[d] & 0\\
0 \ar[r] &\overline{C_{U_{1} \oplus U_{2}}}^{U_{2}}
\ar[rr]\ar[d]&&C_{U_{1}} \ar[rr]\ar[d]&& Q \ar[r]\ar[d] & 0\\
&0 && 0 && 0  & }
$$
where $i_{1}: U_{1} \hookrightarrow U_{1} \oplus U_{2}$ is the natural inclusion.
\medskip\\
Note that the morphism $\langle \psi(g_{1}) \rangle_{U_{1}} \rightarrow \tau^{U_{2}}\big(K_{U_{1} \oplus U_{2}} \big)$ is injective because
\begin{center}
$\begin{array}{lll}
                  \parallel  \langle \psi(g_{1}) \rangle_{U_{1}} \parallel & = & \parallel K_{U_{1}} \parallel
\medskip\\
                   & = & \parallel K_{U_{1} \oplus U_{2}} \parallel, \; \text{by lemma 2.3.2}
\medskip\\
&=& \parallel  \tau^{U_{2}}\big(K_{U_{1} \oplus U_{2}} \big) \parallel, \; \text{because the graded $H^{*}(U_{2})$-module $K_{U_{1} \oplus U_{2}}$ is finite}.
\end{array}$
\end{center}
\medskip
This shows that the graded finite $H^{*}(U_{1} \oplus U_{2})$-module $C_{U_{1} \oplus U_{2}}$ is monogenic generated by an element $\xi \in K_{U_{1} \oplus U_{2}}$. Since the bi-connected Gysin-$(U_{1} \oplus U_{2})$-functor is the restriction of  the Gysin-$V$-functor, ($V = U_{1} \oplus U_{2} \oplus U_{3} $), then by  proposition 2.4.2, we have:
$$\parallel K_{U_{1} \oplus U_{2}} \parallel = \parallel \langle \xi \rangle_{U_{1} \oplus U_{2}} \parallel > \parallel \langle \iota_{1,2} \rangle_{U_{1} \oplus U_{2}} \parallel.$$
This implies an isomorphism of $H^{*}(U_{1} \oplus U_{2})$-modules: $K_{U_{1} \oplus U_{2}} \cong \langle \iota_{1,2} \rangle_{U_{1} \oplus U_{2}} \oplus \langle \xi \rangle_{U_{1} \oplus U_{2}}$.
\medskip\\
By analysing the previous diagram {\bf (D)}, we verify that:
\begin{center}
  \begin{itemize}
    \item [\bf I2.1]  $Im (\overline{j}^{U_{2}}) \cong \overline{\langle \iota_{1,2} \rangle _{U_{1} \oplus U_{2}}}^{U_{2}} \cong \langle \iota_{1} \rangle_{U_{1}} \cong \big(t_{1} \big)_{0}^{n_{1}}\iota_{1}$, $n_{1} \geq 1$, (see $I1.1$).
\medskip
    \item [\bf I2.2]  $\overline{C_{U_{1} \oplus U_{2}}}^{U_{2}} \cong \overline{\langle \xi \rangle_{U_{1} \oplus U_{2}}}^{U_{2}}$ and $ \langle \psi(g_{1}) \rangle_{U_{1}} \cong \tau^{U_{2}}\big(\langle \xi \rangle_{U_{1} \oplus U_{2}}\big)$
\medskip
    \item [\bf I2.3]  $Q \cong \tau^{U_{2}}\big(\langle \iota_{1,2} \rangle_{U_{1} \oplus U_{2}}\big) \cong \big(t_{1} \big)_{0}^{m_{1}}\psi(g_{2})$, $m_{1} \in
\mathbb{N}$, as a graded finite monogenic $H^{*}U_{1}$-module (see notations 2.1).
  \end{itemize}
\end{center}
{\bf I3. The contradiction}.
\smallskip\\
The last line of the previous diagram {\bf (D)}, which is an exact sequence of graded finite monogenic $H^{*}U_{1}$-modules, can now be written, using $I1.2$, as follows:
$$  \xymatrix{0 \ar[r] &\overline{C_{U_{1} \oplus U_{2}}}^{U_{2}} \ar[r] &C_{U_{1}} \ar[d]^-{\cong} \ar[r]^-{\psi_{\mid}} & \tau^{U_{2}}\big(\langle \iota_{1,2} \rangle_{U_{1} \oplus U_{2}}\big) \ar[d]^-{\cong} \ar[r]& 0 \\
& & \big(t_{1} \big)_{0}^{l_{1}}g_{2} & \big(t_{1} \big)_{0}^{m_{1}}\psi(g_{2}) }$$
So we get: $\parallel \overline{C_{U_{1} \oplus U_{2}}}^{U_{2}} \parallel   =  \parallel C_{U_{1}} \parallel > \parallel  \tau^{U_{2}}\big(\langle \iota_{1,2} \rangle_{U_{1} \oplus U_{2}}\big) \parallel$.
\medskip\\
This is equivalent to: $\parallel \big(t_{1} \big)_{0}^{l_{1}}g_{2} \parallel = l_{1} + k_{2} > \parallel \big(t_{1} \big)_{0}^{m_{1}}\psi(g_{2})  \parallel = m_{1} + k_{2}$, $k_{2}$ is the degree of  $g_{2}$. We have then, \mathversion{bold} $l_{1} > m_{1}$  \mathversion{normal}.
\bigskip\\
In conclusion, we have:
$\left\{
  \begin{array}{ll}
    \overline{\langle \iota_{1,2} \rangle _{U_{1} \oplus U_{2}}}^{U_{2}} \cong  \big(t_{1} \big)_{0}^{n_{1}}\iota_{1},\; \text{see}\; I2.1, & \hbox{}
\medskip\\
    \tau^{U_{2}}\big(\langle \iota_{1,2} \rangle_{U_{1} \oplus U_{2}}\big) \cong \big(t_{1} \big)_{0}^{m_{1}}\psi(g_{2}), & \hbox{}
\medskip\\
  m_{1} < n_{1},\; \text{because}\;  m_{1} < l_{1} \leq n_{1},\; \text{see }\; I1.2 & \hbox{}
  \end{array}
\right.$
\bigskip\\
The lemma 2.3.4 (see also the remark 2.3.5) shows the equality of dimensions:
$$ d\big(\overline{\langle \iota_{1,2} \rangle _{U_{1} \oplus U_{2}}}^{U_{2}}\big) =  n_{1} = d\big( \tau^{U_{2}}\big(\langle \iota_{1,2} \rangle_{U_{1} \oplus U_{2}}\big)\big) = m_{1},$$
so the contradiction.


\begin{thebibliography}{MML}

\bigskip

\bibitem[AB]{AB} \textsc{ A. Adem and W. Browder}; The free rang of symmetry of (S$^{n}$)$%
^{k}$; Inventiones Mathematicae 92, 1988, 431-440.

\bibitem[BC]{BC} \textsc{D. J. Benson and J. F. Carlson}; Complexity and multiple complexes, Math. Zeit. 195 (1987), 221-238.

\bibitem[C1]{C1} \textsc{G. Carlsson}; On the non-existence of free action of elementary
abelian groups on products of spheres, American J. of Math. Vol. 102, No
6, 1147-1157 \ (1980).

\bibitem[C2]{C2} \textsc{G. Carlsson}; On the rank of abelian groups acting freely on $(S^{n})^{k}$, Invent. Math. 69 (1982), 393-400.

\bibitem[C3]{C3} \textsc{G. Carlsson}; Free $(\mathbb{Z}/2\mathbb{Z})^{k}$-actions and problems in commutative algebra, Lecture Notes in Math., 1217, Springer, Berlin,(1986), 79-83.

\bibitem[C4]{C4} \textsc{G. Carlsson}; Free $(\mathbb{Z}/2\mathbb{Z})^{3}$-actions on
finite complexes, Algebraic topology and Algebraic K-Theory, W. Browder
(Ed.), Ann. of Math. Studies nº 113, Princeton University Press, Princeton
(1987), 332-344.

\bibitem[Co]{Co} \textsc{P.E. Conner }; On the action of a finite group on $S^{n} \times S^{n}$, Ann. of Math. 66 (1957) 586-588.

\bibitem[DV]{DV} \textsc{J.A. Daccach and J.P. Vieira; }Finite group actions on products
of spheres, manuscripta math. 91, 511-523, (1996).

\bibitem[Hal]{Hal} \textsc{S. Halperin;} Rationnal homotopy and torus actions. Aspects of topology, 293-306, London Math. Soc. Lecture Note Ser., 93 Cambridge Univ. Press, (1985).

\bibitem[Han]{Han} \textsc{B. Hanke;} The stable free rank of symmetry of products of spheres. Invent. Math. 178 (2009), 265-298.

\bibitem[He]{He} \textsc{A. Heller;} A note on spaces with operators. Illinois
J. Math. 3 (1959), 98-100.

\bibitem[MTW]{MTW} \textsc{I. Madsen, C.Thomas and C.T.C.Wall; }The topological space
form problem II, Topology 15 (1976), 375-382.

\bibitem[M]{M} \textsc{J. Milnor; }Groups which act on S$^{n}$ without fixed-points,
Am. J. Math. 79 (1957), 623-630.

\bibitem[MB]{MB} \textsc{J.W. Morgan and H. Bass; }The Smith conjecture, Pure and
applied Mathematics, Academic Press, $n^{\circ}$ 112, 1984.

\bibitem[OY]{OY} \textsc{Osman Berat Okutan and Ergün Yalcin; } Free actions on products of spheres at high dimensions, arXiv: 1207.2363v1, 2012.

\bibitem[R]{R} \textsc{M. Refai; } Group actions on finite
CW-complexes, Indian J. Pure Appl. Math. 24, n° 4, (1993) 245-255.

\bibitem[Sm]{Sm} \textsc{P.A. Smith; }Transformations of finite period,

\qquad \qquad \qquad \qquad \qquad \textbf{I} Ann. of Math.\ 39, (1938),
127-164.

\qquad \qquad \qquad \qquad \qquad \textbf{II} Ann. of Math.\ 40, (1939),
690-711.

\qquad \qquad \qquad \qquad \qquad \textbf{III} Ann. of Math.\ 42, (1941),
446-458.

\bibitem[Sp]{Sp} \textsc{E.H.\ Spanier;} Algebraic Topology, 1966.

\bibitem[TD]{TD}  \textsc{T. tom Dieck,} Transformation groups, De Gruyter Studies in Mathe-
matics 8, 1987.

\bibitem[Z]{Z} \textsc{S. Zarati; }D\'{e}faut de stabilit\'{e} d'op\'{e}rations
cohomologiques, Th\`{e}se de 3-\`{e}me cycle, Publications d'Orsay $n^{\circ}$ 78-07.

\end{thebibliography}
\end{document}